\documentclass{article}

\usepackage{amsmath}
\usepackage{graphicx}
\usepackage{url}

\usepackage{amsmath, amssymb, amstext, amsopn, amsxtra, graphicx, color, calligra, hyperref, url, ulem, polynom, permute, ifthen, float, epstopdf, linsys, marvosym, rotating}

\newtheorem{theorem}{Theorem}
\newtheorem{corollary}[theorem]{Corollary}

\begin{document}

\newcommand{\mob}{M\"obius}
\newcommand{\murraymod}[3]{{#1 \equiv #2 \pmod{#3}}}  
\newcommand{\pf}{{\bf Proof}:  }
\newcommand{\qed}{\mbox{\(\square\,\,\)}}
\newcommand{\field}[1]{\ensuremath{\mathbb #1}}
\newcommand{\CC}{\field C}
\newcommand{\FF}{\field F}
\newcommand{\NN}{\field N}
\newcommand{\ZZ}{\field Z}
 
\title{\mob\ Polynomials}

\author{Will Murray\footnote{I thank Art Benjamin and Robert Mena for suggesting this investigation and many useful subsequent thoughts. The referees also offered many helpful ideas that dramatically improved this paper.}\\
California State University, Long Beach\\
Long Beach, CA 98040-1001\\
Will.Murray@csulb.edu}
\date{December 11, 2013}

\maketitle

What do bracelets, juggling patterns, and irreducible polynomials over finite fields have in common? We can count all of them using the same techniques. In this paper, we identify a set of polynomials that unify all the counting problems above. Despite their wide application, they have not received an explicit name in the literature. We will christen them \textit{\mob\ polynomials} and explore some of their properties. We will highlight the role that \mob\ polynomials have played in the contexts above and then use them to 
extend a classic combinatorial proof from Fermat's little theorem to
Euler's totient theorem.

In the first section below, we will define our polynomials and derive some key facts about them. After a brief digression to enjoy the graphs of our polynomials in the complex plane, we will see that $M_n(x)$ gives the number of aperiodic bracelets of length $n$ that can be built using $x$ possible types of gems. An immediate corollary will be that \mbox{$\murraymod{M_n(x)}0n$} for all $x \in \ZZ$. In three subsequent sections, we will apply our polynomials to count juggling patterns, to count irreducible polynomials over finite fields, and to prove Euler's totient theorem.

\section{Definition and properties}

To construct our polynomial, we first recall the \textit{\mob\ $\mu$} function defined on a positive integer $n = p_1^{e_1}\cdots p_r^{e_r}$, where the $p_i$ are distinct primes:
\begin{eqnarray*}
\mu(1) & := & 1\\
\mu\left( n = p_1^{e_1}\cdots p_r^{e_r} \right) & := & \begin{cases}(-1)^r & \text{if all } e_i = 1\\ 0 & \text{if any } e_i > 1\end{cases}
\end{eqnarray*}
Number theorists use $\mu$ for \textit{\mob\ inversion}, as in \cite{murray:burton}:
If $f$ and $g$ are functions such that 
\[
\sum_{d|n}f(d) = g(n),
\]
then we can solve for $f$ in terms of $g$ via $\mu$:
\[
f(n) = \sum_{d|n}\mu\left( \frac nd \right) g(d).
\]
We can now meet our main object of study.

\smallskip
\textsc{Definition.\ } For each integer $n\geq 1$, the $n$th \textit{\mob\ polynomial} is defined to be
\[
M_n(x) := \sum_{d|n} \mu\left( \frac nd \right) x^d.
\]
\smallskip

\noindent For example, when $n=12$, we get $M_{12}(x) = x^{12}-x^6-x^4+x^2$.

The first nice property of the \mob\ polynomial is that its coefficients are all 0 or $\pm 1$. The coefficient of $x^d$ in $M_n(x)$ is 0 if any square divides $\frac nd$; therefore the only nonzero terms that appear in $M_n(x)$ are those corresponding to divisors $d|n$ for which $\left.p_1^{e_1-1}\cdots p_r^{e_r-1}\right|d$, where $n = p_1^{e_1}\cdots p_r^{e_r}$. Thus the leading term of $M_n(x)$ is always $x^n$ and the smallest nonzero term is $(-1)^r x^d$, where $d = p_1^{e_1-1}\cdots p_r^{e_r-1}$. The multiplicity of the root at 0 is therefore $p_1^{e_1-1}\cdots p_r^{e_r-1}$.

We note also that $M_n(x)$ always has $2^r$ nonzero terms, corresponding to the $2^r$ divisors between $p_1^{e_1-1}\cdots p_r^{e_r-1}$ and $n$, one for each subset of $\{p_1,\dots,p_r\}$. If $n>1$, then for exactly half of these divisors, $\frac nd$ factors into an odd number of primes, so exactly half of the coefficients are $-1$ and the other half are $1$. This proves that if $n>1$, then $M_n(1) = 0$.

This is the beginning of an interesting line of study. It turns out that many \mob\ polynomials have zeroes at many roots of unity. Here we will just investigate $M_n(-1)$. For $n = 1$ and $n=2$, we have $M_1(x) = x$ and \mbox{$M_2(x) = x^2 - x$}, so $-1$ is not a zero, but it is a zero for all higher \mob\ polynomials:

\begin{theorem}
If $n > 2$, then $M_n(-1) = 0$.
\end{theorem}

We present two proofs.

\pf  [Straightforward proof:] We will examine each term of $M_n(-1)$ and see that exactly half are negative. There are two cases:

First suppose that $4|n$. As we saw above, the only divisors $d|n$ that appear are those for which $p_1^{e_1-1}\cdots p_r^{e_r-1}|d$. But all of these divisors are even, so $M_n(x)$ contains only even powers of $x$. Thus, $M_n(-x) = M_n(x)$, and in particular $M_n(-1) = M_n(1) = 0$.

Next suppose that $4\nmid n$. Then since $n > 2$, there must be at least one odd prime $p$ dividing $n$. Then among the $d$ for which $\left.p_1^{e_1-1}\cdots p_r^{e_r-1}\right|d$, half contain the final power of that prime $p$ as a factor and half do not. We can pair each $d$ that does not contain that final $p$ with $pd$ and note that $\mu\left( \frac nd \right) (-1)^d$ cancels with $\mu\left( \frac n{pd} \right) (-1)^{pd}$, since $\mu\left( \frac nd \right) $ and $\mu\left( \frac n{pd} \right) $ have opposite signs but $(-1)^d$ and $(-1)^{pd}$ have the same sign. For example, in $M_6(x) = x^6 - x^3 - x^2 + x$, we cancel $x^6$ with $-x^2$ and $-x^3$ with $x$. Therefore the total sum is 0. 
\qed

\pf  [Analytic proof:] We first note that since $r$ is the number of \textit{distinct} primes dividing $n$, we have $n > 2^r$. (The only case when this inequality would not be strict would be if $r = 1$ and $n = 2$, but we excluded that case by hypothesis.) Now, $M_n(-1)$ has $2^r$ terms, each of which is $\pm 1$, so $\left|M_n(-1)\right| \leq 2^r < n$. But we will see in the corollary to our theorem on bracelets, which we will prove independently in a later section, that $\murraymod{M_n(x)}0n$ for all $x \in \NN$, so $\murraymod{M_n(-1) \equiv M_n(n-1)}0n$. Thus, the only possible value of $M_n(-1)$ is 0.
\qed

We will illustrate these results in the next section before moving on to combinatorial applications.

\section{Digression: The graphs of \mob\ polynomials}

As mentioned above, many \mob\ polynomials have zeroes at many roots of unity.
It is interesting to consider a \mob\ polynomial as a function of a complex argument $z$ and to examine its values when $z$ lies on the unit circle in $\CC$. We can graph them by evaluating $M_n(z)$ for $z = e^{i\theta}$, \mbox{$0 \leq \theta \leq 2\pi$}, and plotting the results in the complex plane. For example, 
Figure~1 shows the graphs of \mbox{$M_{15}(z) = z^{15} - z^5 - z^3 + z$} and \mbox{$M_{17}(z) = z^{17} - z$}.
\begin{figure}[tpb] 
\begin{center}
\resizebox!{60mm}{\includegraphics{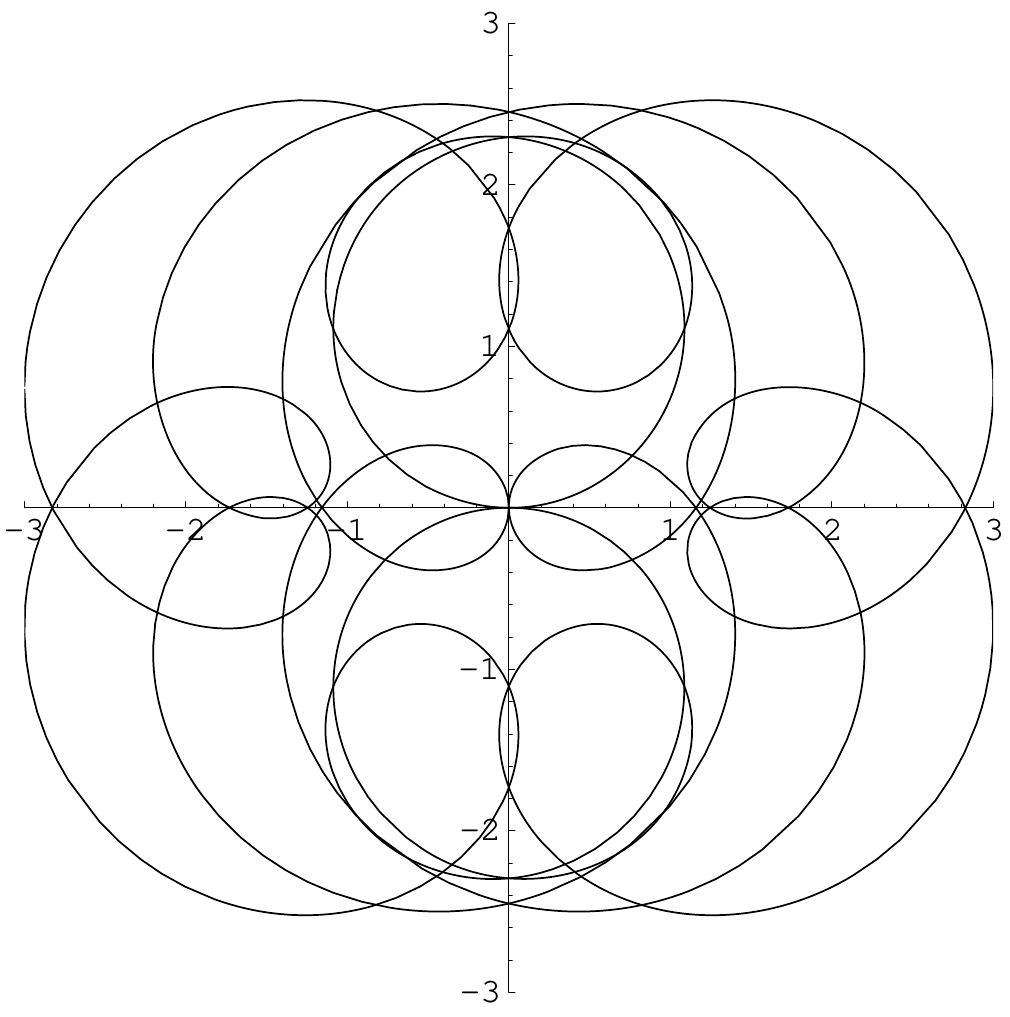}}
\resizebox!{60mm}{\includegraphics{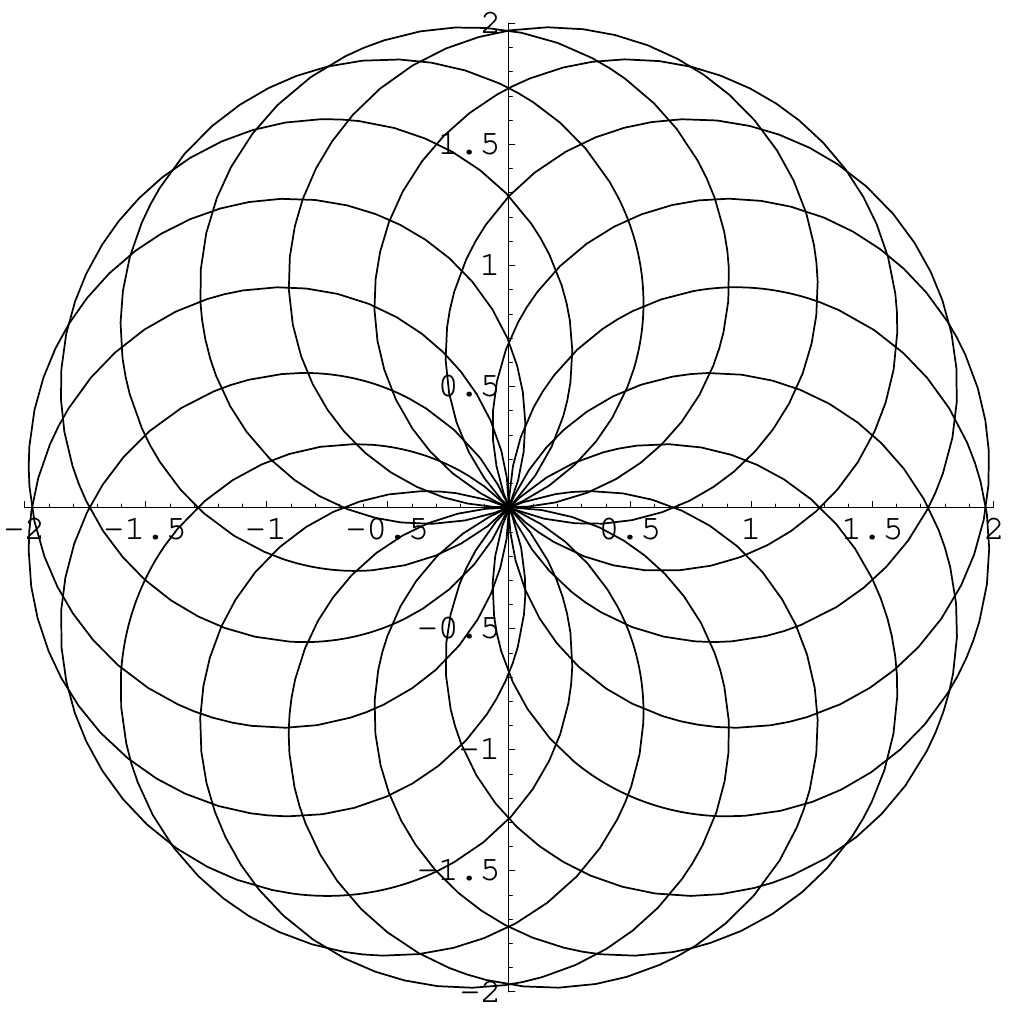}}
\end{center}
\caption{The graphs of $M_{15}(z)$ and $M_{17}(z)$ for $z = e^{i\theta}$, $0 \leq \theta \leq 2\pi$.}
\label{murrayfig-spirograph}
\end{figure}

The graphs cross the origin at zeroes of the polynomials. $M_{15}(z)$ has zeroes at $\pm 1$ and $\pm i$, while $M_{17}(z) = z^{17} - z = z\left(z^{16}-1\right)$ has zeroes at all sixteenth roots of unity.

The symmetries of the graphs reflect the structure of the polynomials. $M_{15}(z)$ satisfies $M_{15}(-z)= -M_{15}(z)$, giving the graph rotational symmetry around the origin; $M_{15}(\bar z)=\overline{M_{15}(z)}$, giving the graph vertical symmetry across the real axis; and $M_{15}(-\bar z)=-\overline{M_{15}(z)}$, giving the graph horizontal symmetry across the imaginary axis. $M_{17}(z)$, on the other hand, satisfies 
\[
M_{17}(\omega z) = \omega z\left( (\omega z)^{16}-1 \right) = \omega z\left( z^{16}-1 \right) = \omega M_{17}(z),
\]
for the sixteenth root of unity $\omega = e^{\frac{2\pi i}8}$, so the graph is symmetric with respect to rotation through $\omega$.

There are many more beautiful symmetries and patterns to be discovered by graphing \mob\ polynomials and, more generally, other functions $f:\CC \rightarrow \CC$ on the unit circle. I invite you to play with them yourself. In the meantime, we will return to combinatorial applications of \mob\ polynomials. 

\section{\mob\ polynomials and bracelets}
We are ready for our first combinatorial result on \mob\ polynomials. We would like to build circular bracelets of length $n$ using $x$ possible types of gems. We can think of each bracelet as a word of $n$ letters, and we have $x$ choices for each letter, so there are $x^n$ possible bracelets in all. However, we wish to exclude those bracelets that are \textit{periodic} with respect to any proper divisor $d|n$, that is, those that after rotating by $d$ gems look the same as themselves. For example, the bracelets XOXOXO and OXOOXO are periodic with periods 2 and 3, respectively, but the bracelet XOXOOO is aperiodic.

\begin{theorem}
The \mob\ polynomial $M_n(x)$ gives the number of aperiodic bracelets of length $n$ with $x$ possible types of gems.
\end{theorem}

\pf  Note that every bracelet of length $n$ is periodic with respect to some \textit{fundamental} (that is, shortest) period $d|n$; the aperiodic bracelets are those for which $d=n$. For each $d|n$, there is a one to one correspondence between the bracelets of length $n$ with fundamental period $d$ and the aperiodic bracelets of length $d$. Let us define $g(d)$ to be the number of aperiodic bracelets of length $d$. Then our correspondence gives us $\sum_{d|n} g(d) = x^d$. We now apply \mob\ inversion:
\[
g(n) = \sum_{d|n} \mu\left( \frac nd \right) x^d = M_n(x)
\]
This proves the theorem.
\qed

\begin{corollary}
The \mob\ polynomial satisfies $\murraymod{M_n(x)}0n$ for all $x \in \NN$.
\end{corollary}

\pf  We can sort the aperiodic bracelets of length $n$ into groups that look like each other after rotation. For example, with $n=6$ and $x=2$, one group would be $\{$XOXOOO, OXOXOO, OOXOXO, OOOXOX, XOOOXO, OXOOOX$\}$. Because we are only considering aperiodic bracelets, each group contains exactly $n$ bracelets, proving the corollary.
\qed

We will recycle this argument for several applications in the following sections, first to count juggling patterns and irreducible polynomials and then to prove Euler's totient theorem. We note here that Bender and Goldman \cite{murray:bender} derive a similar formula in which they count the total number of rotationally distinct bracelets of length $n$ (including periodic ones), meaning for example that they consider OXOOXO and OOXOOX to be the same. Thus each divisor $d|n$ contributes
\[
\frac{1}{d} M_d(x) = \frac{1}{d} \sum_{c|d}\mu\left( \frac dc \right) x^c
\]
bracelets, giving a total of 
\[
\sum_{d|n} \frac{1}{d} \sum_{c|d}\mu\left( \frac dc \right) x^c
\]
bracelets in all. This expression simplifies pleasantly when we use Gauss's identity $n = \sum_{d|n} \phi(n)$, giving
\[
\frac{1}{n} \sum_{d|n} \phi\left( \frac nd \right) x^d.
\]
We omit the details because we will not need them here.

\section{\mob\ polynomials and juggling patterns}
Jugglers and mathematicians describe juggling patterns using \textit{siteswap} notation, in which strings of nonnegative integers represent throws to different heights. For example, 441 represents throwing two balls high and then quickly passing a third ball from one hand to the other underneath them. (Actually, the hand order, and even the number of hands, are irrelevant to the notation. Each positive integer just represents the number of time beats from when a ball is thrown to when it is thrown again, and a zero represents a beat in which no ball is thrown.) Jugglers usually repeat a sequence of throws periodically, so 441 is shorthand for a pattern in which one would throw red, blue, and green balls in the order R G B B R G G B R R G B B R G G B R \dots.

There are many websites that animate juggling patterns, notably Boyce's Juggling Lab~\cite{murray:boyce}. Siteswap notation has spread like wildfire from academia to the mainstream juggling community of hobbyists, performers, and competitors. In 2005 for example, Japanese street performer Kazuhiro Shindo won the International Jugglers' Association championships with a routine based on variations on 7441, and he even shaved the formula into the back of his head (Figure~2).

\begin{figure}[htbp]
\begin{center}
\resizebox{40mm}!{\includegraphics[trim = 7cm 3cm 6.5cm 3cm, clip]{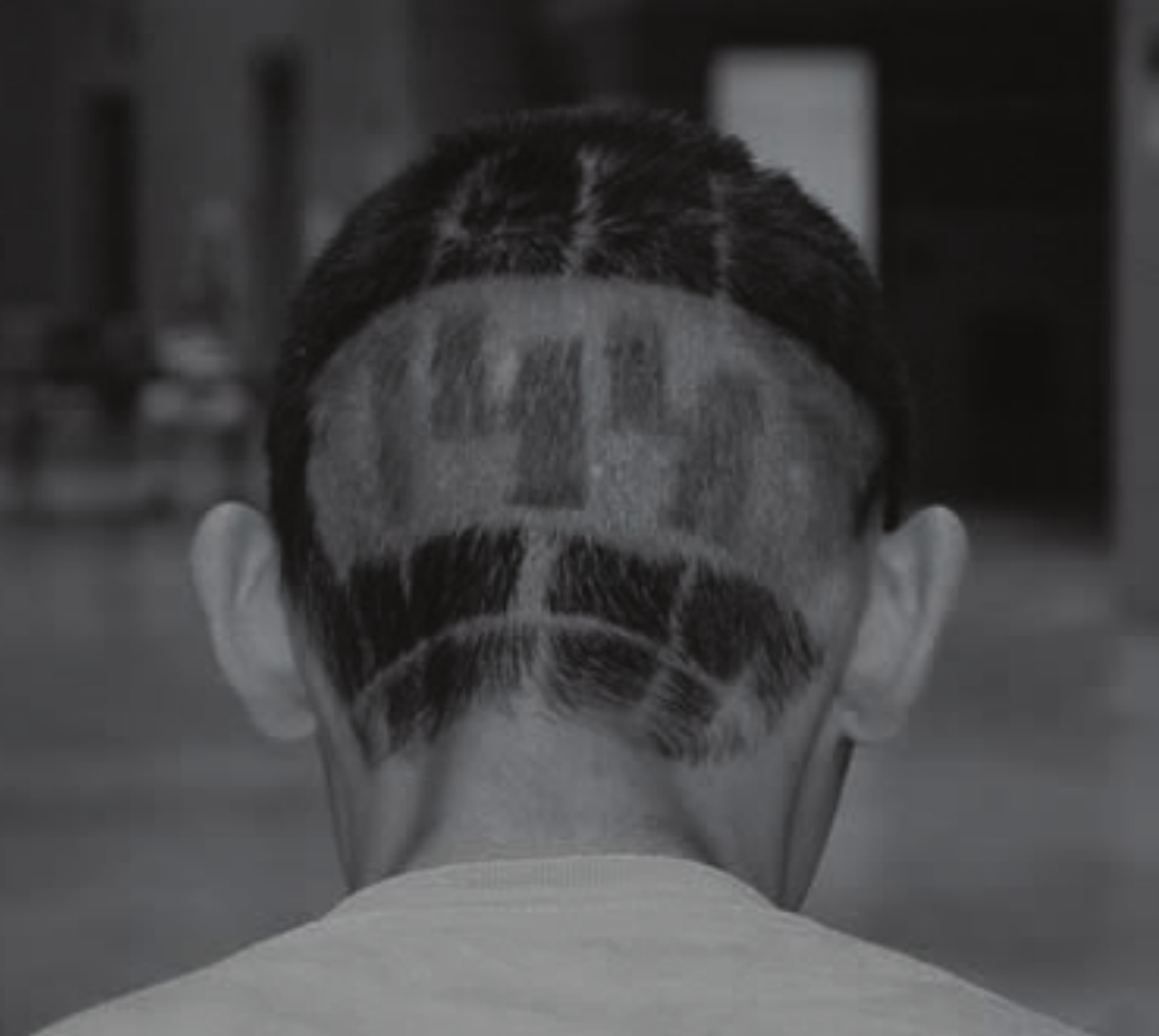}} 
\end{center}
\caption{Kazuhiro Shindo and 7441 (photo courtesy of Joyce Howard)}
\label{murrayfig-shindo}
\end{figure}

Academically minded jugglers have long known that a string of nonnegative integers $a_1\cdots a_n$ is a valid siteswap pattern if and only if for all $1 \leq i \neq j \leq n$, we have $i + a_i \not\equiv j + a_j$ (mod $n$). This condition ensures that the balls do not collide upon landing, since throw $i$ lands at time $i + a_i$ (mod $n$). In their seminal paper \cite{murray:buhler}, Buhler, Eisenbud, Graham, and Wright prove (via some nontrivial combinatorics) the remarkable theorem that the number of patterns of period $n$ with strictly fewer than $b$ balls is exactly~$b^n$.

As an example, the predicted $2^4 = 16$ patterns of period four with zero or one ball(s) are 0000, 4000, 0400, 0040, 0004, 3001, 1300, 0130, 0013, 2020, 0202, 1111, 2011, 1201, 1120, and 0112. Actual jugglers, however, would only list 4000, 3001, and 2011, since all the others are either cyclic copies of these three (recall that a juggler repeats a pattern indefinitely without caring where it begins and ends) or have fundamental period less than four. For example, 2020 has fundamental period two.

A juggler would define $f(n,b)$ to be the number of siteswap patterns of fundamental period $n$, where cyclic copies of patterns such as 3001 and 0130 are considered the same. Then for each divisor $d|n$, the total of $b^n$ includes patterns of fundamental period $d$, and each pattern is counted $d$ times because it can be rotated $d$ places until it repeats itself. For example, the period two pattern 20 is counted twice in the list of period four patterns above. The total thus breaks down as follows:
\[
\sum_{d|n} df(d,b) = b^n
\]
Let us define the temporary function $T(d,b) := df(d,b)$. Then our formula becomes $\sum_{d|n} T(d,b) = b^n$, and \mob\ inversion gives us 
\[
T(n,b) = \sum_{d|n}\mu\left( \frac nd \right) b^d = M_n(b),
\]
that is,
\[
f(n,b) = \frac{1}{n} M_n(b).
\]
For example, the count of patterns of period four with zero or one ball(s) is
\begin{eqnarray*}
\frac{1}{4} M_4(2) 
 &=& \frac{1}{4} \sum_{d|4} \mu\left( \frac4d \right) 2^d \\
 &=& \frac{1}{4}\left[\underbrace{\mu(4) 2^1}_{d=1} + \underbrace{\mu(2) 2^2}_{d=2} + \underbrace{\mu(1) 2^4}_{d=4}\right] \\
 &=& \frac{1}{4} \left[ 0 - 4 + 16 \right] = 3
\end{eqnarray*}
confirming the three distinct patterns 4000, 3001, and 2011.

The original four authors, all accomplished jugglers themselves, gave this formula in \cite{murray:buhler}. In fact, they derived the formula for the number of patterns with \textit{exactly} $b$ balls:
\begin{eqnarray*}
f(n,b+1) - f(n,b) 
 &=& \frac{1}{n} M_n(b+1) - \frac{1}{n} M_n(b) \\
 &=& \frac{1}{n} \sum_{d|n} \mu\left(\frac nd\right) \left[ (b+1)^d - b^d \right]
\end{eqnarray*}
(This expression is denoted $M(n,b)$ in \cite{murray:buhler}; we will avoid this notation because it could be confused with our $M_n(b)$.)

For example, the number of period-three patterns with exactly three balls is
\[
\frac{1}{3}\left[\underbrace{\mu(3) \left( 4^1-3^1 \right) }_{d=1} + \underbrace{\mu(1) \left( 4^3-3^3 \right) }_{d=3} \right] = 12,
\]
and the patterns are 423, 441, 450, 522, 531, 603, 612, 630, 711, 720, 801, and 900.

\section{\mob\ polynomials and irreducible polynomials}

\mob\ polynomials find another application in counting irreducible polynomials over finite fields. In fact, with a little finite field theory, we can give a formula using the same proof as for juggling patterns:

\begin{theorem}
The number of monic irreducible polynomials of degree exactly $n$ over the finite field $\FF_p$ is
\[
\frac{1}{n} M_n(p) = \frac{1}{n} \sum_{d|n}\mu\left( \frac nd \right) p^d.
\]
\end{theorem}

\pf  Define $f(n)$ to be the number of such polynomials. Each element \mbox{$\alpha \in \FF_{p^n}$} satisfies an irreducible monic polynomial of some degree $d$ over $\FF_p$, and we know $d|n$ because $\FF_p(\alpha) = \FF_{p^d} \subseteq \FF_{p^n}$. In this list of irreducible polynomials, each irreducible polynomial of degree $d$ is counted $d$ times, once for each of its roots. Therefore, $\sum_{d|n} df(d) = p^n$, and just as in the derivation of the juggling formula, we get $f(n) = \frac{1}{n} \sum_{d|n}\mu\left( \frac nd \right) p^d = \frac{1}{n} M_n(p)$.
\qed

In fact, the same proof gives a more general version of the theorem, that the number of monic irreducible polynomials of degree exactly $n$ over the finite field $\FF_{p^e}$ is
\[
\frac{1}{n} M_n\left( p^e \right) = \frac{1}{n} \sum_{d|n}\mu\left( \frac nd \right) p^{de}.
\]

This formula, of course, is known in the literature, for example in Dornhoff and Hohn \cite{murray:dornhoff}. In \cite{murray:chebolu}, Chebolu and Min\'a\v{c} give a proof based on inclusion-exclusion that interprets each term in the sum above in terms of field theory.

\section{\mob\ polynomials and Euler's totient theorem}

Our corollary on bracelets above gives us an immediate combinatorial proof of Fermat's little theorem.  Take $n=p$ to be prime, and we get $\murraymod{M_p(x) = x^p - x}0p$.  This proof of Fermat appeared first in Dickson \cite{murray:dickson} and later with many interesting variations in Anderson, Benjamin, and Rouse \cite{murray:anderson}.  In their master compendium of combinatorial proofs \cite{murray:benjamin}, Benjamin and Quinn issue the challenge:  ``Although we do not know of a combinatorial proof of [Euler's totient theorem that $\murraymod{a^{\phi(n)}}1n$ when $a$ and $n$ are relatively prime], we would love to see one!''  
We can use our corollary to prove a special case of Euler, when $n=p^e$.  In this case, Euler's theorem becomes $\murraymod{a^{p^e - p^{e-1}}}1{p^e}$, or $\murraymod{a^{p^e}}{a^{p^{e-1}}}{p^e}$.

\pf  [Proof of Euler's totient theorem (special case):] We evaluate the \mob\ polynomial $M_{p^e}(x)$ at $x=a$:
\[
M_{p^e}(a) := \sum_{d|n} \mu\left(\frac {p^e}d \right) a^d
\]
This polynomial has only two nonzero terms, which are $a^{p^e} - a^{p^{e-1}}$. By the corollary on bracelets, we have $\left. p^e \right| \left( a^{p^e} - a^{p^{e-1}} \right)$, proving the special case.
\qed

It is easy to derive the general version of Euler from the combinatorial special case as follows: Suppose $a$ and $n = p_1^{e_1}\cdots p_r^{e_r}$ are relatively prime, where the $p_i$ are distinct primes. Then for each prime $p_i$, we have:
\begin{align*}
&\left. p_i^{e_i} \right| \left( a^{p_i^{e_i}} - a^{p_i^{e_i-1}} \right) && \text{by the combinatorial special case}\\
&\left. p_i^{e_i} \right| \left[ a^{p_i^{e_i-1}} \left( a^{p_i^{e_i}-p_i^{e_i-1}}-1 \right) \right] && \text{by factoring}\\
&\left. p_i^{e_i} \right| \left( a^{p_i^e-p_i^{e_i-1}}-1 \right) && \text{since $\left( a,p_i \right) = 1$}\\
&\left. p_i^{e_i} \right| \left( a^{\phi\left(p_i^{e_i}\right)}-1 \right) &&
\end{align*}
We know that $\phi$ is multiplicative on relatively prime arguments, so \linebreak\mbox{$\phi(n) = \phi\left(p_1^{e_1}\right) \cdots \phi\left(p_r^{e_r}\right)$}. In particular, each $\left.\phi\left(p_i^{e_i}\right)\right|\phi(n)$, so 
\[
\left. \left( a^{\phi\left(p_i^{e_i}\right)}-1 \right) \right| \left( a^{\phi(n)}-1 \right),
\]
and transitivity yields \mbox{$p_i^{e_i} \left|\left( a^{\phi(n)}-1 \right) \right.$}. Combining these gives us $n\left|\left(a^{\phi(n)}-1\right)\right.$, the general version of Euler's theorem.

{\bf Summary}:  We introduce the \mob\ polynomial $ M_n(x) = \sum_{d|n} \mu\left( \frac nd \right) x^d $, which gives the number of aperiodic bracelets of length $n$ with $x$ possible types of gems, and therefore satisfies $\murraymod{M_n(x)}0n$ for all $x \in \ZZ$. We derive some key properties, analyze graphs in the complex plane, and then apply \mob\ polynomials combinatorially to juggling patterns, irreducible polynomials over finite fields, and Euler's totient theorem.
\label{Murraypage:last}

\vfill\eject
\end{document}